\documentstyle[12pt]{article}

\newtheorem{thm}{Theorem}
\newtheorem{prop}{Proposition}

\newtheorem{lemma}{Lemma}

\newcommand{\QED}{\hfill $\Box$}

\newcommand{\R}{\mbox{{\bf{R}}}}
\newcommand{\Z}{\mbox{{\bf{Z}}}}
\newcommand{\ext}{\bigwedge\hspace{-.02in}}

\title{Local Properties of Self-Dual Harmonic 2-forms on a
4-Manifold}
\author{Ko Honda}
\date{May 27, 1997}

\begin{document}
\maketitle

\begin{abstract}
We will prove a Moser-type theorem for self-dual harmonic 2-forms on closed
4-manifolds, and use it to classify local forms on neighborhoods of 
singular circles on which the 2-form vanishes.
Removing neighborhoods of the circles, we obtain a symplectic manifold
with contact boundary - we show that the contact form on each $S^1\times
S^2$, after a slight modification,  must be one of two possibilities.
\end{abstract}

\section{Introduction}
This paper is a study of generic self-dual (SD) harmonic 2-forms $\omega$
near its zero set. 
Let $M^4$ be a closed, oriented 4-manifold with 
$b^+_2>0$.  Then, for a pair $(\omega,g)$ consisting of a generic
metric $g$ and a self-dual harmonic 2-form $\omega$ with respect to
$g$, $(\omega,g)$ represents a section of $\bigwedge_g^+\rightarrow
M$, which is transverse to the zero section.  Here $\bigwedge_g^+$
is the subbundle of $\bigwedge^2TM\rightarrow M$ whose fiber over a
point $p\in M$ is $\bigwedge^+_g(p)=\{\omega|*_g\omega=\omega\}$.
In particular, the zeros of $\omega$ are disjoint embedded circles.
Since $\omega\wedge\omega=\omega\wedge*\omega$, $\omega$ is
nondegenerate at $p$ if and only if $\omega(p)\not = 0$.  That is,
$\omega$ is closed and symplectic away from the union of circles
$C$, and is identically $0$ on $C$.  

We also have a relative version of the previous discussion, which is the 
following theorem (cf. \cite{Honda}):

\begin{thm}
Let $(\omega_0,g_0)$ and $(\omega_1, g_1)$ be  harmonic
forms transverse to $\ext^+_{g_0}$ and $\ext^+_{g_1}$, respectively.  If there
exists a path $(\omega_t, g_t)$ of harmonic forms $\omega_t$ with
respect to $g_t$ such that $\omega_t\not=0$ for all $t\in[0,1]$,
then there exists a $G_\delta$-set of perturbations
$\{(\widetilde\omega_t, \tilde{g_t})\}$ of this path, fixing endpoints,
such that $\{(\widetilde\omega_t, \tilde{g_t})\}$ has regular zeros in
$M\times [0,1]$.  
\end{thm}

\noindent
{\bf Note:} The conditions for the theorem are minimal.  The space
$\{(\omega,g)| g\in \mbox{Met}^k(M),$ $ *\omega=\omega,$ $ \Delta_g\omega=0\}$ is
diffeomorphic to $\R^{b_2^+}\times \mbox{Met}^k(M)$, 
where $\mbox{Met}^k(M)$ is the space
of $C^k$-metrics on $M$, and, as long as $b_2^+>1$, we can always
connect $(\omega_0, g_0)$ to $(\omega_1,g_1)$ via a cobordism  such
that $\omega_t\not = 0$ for all $t\in[0,1]$.  In the case $b^+_2=1$,
as long as $(\omega_0, g_0)$ and $(\omega_1,g_1)$ lie on the same
side of the real line, there exists a cobordism.

\vskip.2in
We briefly outline the contents of the paper.
In Section 3, we will discuss a version of Moser's theorem (Theorem
\ref{thm:sddiffeo}) which applies to our {\em singular symplectic forms}.
In Section 4, we classify local normal forms of the singular symplectic
forms near an $S^1$, with an eye towards global results, and in the last 
section we discuss the induced contact structures on the boundaries of
$N(S^1)$.  These remarks lay the groundwork for the 
Floer homology of 
singular symplectic forms, which we hope to return to in a subsequent 
paper.

\section{Almost complex structures}

Observe first that we can define an almost complex structure $J$ on
$M-C$.  

\begin{prop}
If $\omega$ is a self-dual harmonic 2-form which is nondegenerate on
a connected set $M-C$, then there exists a unique almost complex
structure $J$ compatible with $\omega$ and $\tilde{g}$ on $M-C$,
where $\tilde{g}$ is conformally equivalent to $g$.
\end{prop}
\noindent
{\bf Proof:} Any 2-form $\omega$ can be written, with
respect to the metric $g$, as
$$ \omega=\lambda_1e_1e_2+\lambda_2e_3e_4,$$
with $e_1,...,e_4$ orthonormal and positively oriented at a point
$p\in M-C$.  

For $\omega$ to be self-dual, $\lambda_1=\lambda_2$.  Hence,
$$\omega=\lambda(e_1e_2+e_3e_4).$$
This $\lambda$ is well-defined up to sign:  Simply consider ${1\over
2}\omega\wedge\omega=\lambda^2e_1...e_4=\lambda^2dv_g$, with $dv_g$ the 
volume form for $g$.
Since $\lambda^2$ is only dependent on $\omega$ and $g$, we can
determine $\lambda$ up to sign.  However, taking advantage of $M-C$
being connected, we may fix $\lambda$ on all of $M-S$ so that
$\lambda>0$.  

We then set $J:e_1\mapsto e_2, e_2\mapsto -e_1, e_3\mapsto e_4,
e_4\mapsto -e_3.$ 
This definition is equivalent to the following: Let
$\tilde{g}=\lambda g$, and define $J$ such that
$\tilde{g}(x,y)=\omega(Jx,y)$.  Hence we see that if there is a $J$
compatible with $\omega$ and $\tilde{g}$, it must be unique.
Thus $J$ is compatible with $\omega$ and $\tilde{g}=\lambda g$ on
$M-C$.  \QED

\vskip.20in
Observe that $\omega$ is defined on all of $M$ and is zero on $C$,
$\tilde{g}$ can be defined on all of $M$ and is zero on $C$, but is
not smooth on $C$, while $J$ is defined only on $M-C$.

Let $\{(\omega_t,g_t)\}$ be a regular cobordism.  As in the previous
proposition, we can define $\lambda_t$ uniformly over
$\bigcup_{t\in[0,1]}(M\times \{t\}-C_t)$ and get a family
$\{(\omega_t,\tilde{g_t},J_t)\}$, which is compatible where defined.

\vskip.20in

\section{Moser argument for self-dual harmonic 2-forms}

Consider $M^4$ as above.  Let $\{\omega_t\}$ be a generic family of
self-dual harmonic 2-forms such that

\vskip.15in
(i) $[\omega_t]\in H^2(M;\R)$ is constant.

\vskip.15in
(ii) The sets $C_t=\{x\in M|\omega_t(x)=0\}$ are all $S^1$'s; hence via a
diffeomorphism, we may assume that $C=C_t$ is a fixed $S^1$.  

\vskip.15in
(iii) $[\omega_t]\in H^2(M,C;\R)$ does not vary with $t$.

\vskip.15in
\noindent
If we assume that $C$ is contractible, then we are asking for the following:

\vskip.15in
(iii$'$) Let $\Omega$ be an oriented surface with $\partial \Omega=C$.
 Then $\int_\Omega
\omega_t$ does not vary with $t$.

\vskip.15in
\noindent
Then we have the following:
\begin{thm}\label{thm:sddiffeo}
There exists a 1-parameter family of $C^1$-diffeomorphisms of $M$,
which is smooth away from $C$, and takes $(M-C, \omega_0)
{\stackrel\sim\rightarrow}(M-C,\omega_1)$ symplectically.
\end{thm}

This generalizes the classical
\begin{thm}[Moser] Let $\{\omega_t\}$ be a family of symplectic forms on
a closed manifold $M$.  Provided $[\omega_t]\in H^2(M;\R)$ is fixed,
there is a 1-parameter family of diffeomorphisms $\phi_t$ such that
$\phi_t^*\omega_t=\omega_0$.  
\end{thm}

\noindent
{\bf Proof: (Moser)} Let $\eta_t$ be a 1-parameter family of 1-forms
such that ${d\omega_t \over dt}=d\eta_t$. Thus, if we define $X_t$ such
that $i_{X_t}\omega_t=\eta_t$, then ${\cal L}_{X_t}\omega_t=(i_{X_t}\circ
d+d\circ i_{X_t})\omega= d\eta_t$, which, integrated,  gives a
1-parameter family $\phi_t$ such that $\phi_t^*\omega_t=\omega_0$.\QED

\vskip.15in
\noindent
{\bf Proof: (Theorem~\ref{thm:sddiffeo})} The point here is to find a suitable $\eta_t$
such that ${d\omega_t\over dt}=d\eta_t$ and $\eta_t|_C=0$.  Fix some
$\tilde\eta_t$ such that ${d\omega_t\over dt}=d\tilde\eta_t$.  We shall
find a function $f_t$ on $M$ such that $\tilde\eta_t=df_t$ ``up to
first order'' near $C$.  

Condition (iii) implies that there exists an $f_t$ on $C$ such that
$i^*\tilde\eta_t=df_t$, where $i:C\rightarrow M$ is the inclusion,
i.e. $i^*\tilde\eta_t$ is exact.  In particular, assuming (iii$'$) we have
$$\int_Ci^*\tilde\eta_t=\int_{\partial\Omega}i^*\tilde\eta_t=
\int_\Omega d\tilde\eta_t=\int_\Omega{d\omega_t\over dt}=0.$$
In order to extend $f_t$ to a neighborhood $N(C)$ of $C$, first
observe that there is only one orientable rank 3 bundle over $S^1$
($\pi_1(BSO(3))=0$ implies $S^1\rightarrow BSO(3)$ is homotopically
trivial) and hence $N(C)\simeq C\times D^3$.  Choose coordinates
$(\theta, x_1,x_2,x_3)$ such that $d\theta,dx_1,dx_2,dx_3$ at
$(\theta,0)$ are orthonormal.  

Setting 
$$f_t(\theta,x_1,x_2,x_3)=f_t(\theta,0)+\sum_i
\tilde\eta_i(\theta,0)x_i + {1\over 2}\sum_{i,j}{\partial
\tilde\eta_i \over \partial x_j}(\theta,0)x_ix_j$$
on $N(C)$, where $\tilde\eta_t=
\tilde\eta_\theta d\theta+\sum_i\tilde\eta_idx_i$, we have
\begin{eqnarray*} df_t(\theta,x_1,x_2,x_3)&=&{\partial f_t\over
\partial\theta}(\theta,0)d\theta+\sum_i{\partial \tilde\eta_i\over
\partial\theta}(\theta,0)x_id\theta\\
&+&\sum_i\tilde\eta_i(\theta,0)dx_i +{1\over
2}\sum_{i,j}{\partial\tilde\eta_i\over \partial x_j}(\theta,0)
(x_idx_j+x_jdx_i)
\end{eqnarray*}
up to first order in the $x_i$'s.  Now observing that 
\begin{equation} {\partial f\over
\partial\theta}(\theta,0)=\tilde\eta_\theta(\theta,0),
\end{equation}
\begin{equation} \label{eq2}d\tilde\eta_t(\theta,0)=0,
\end{equation}
and that Equation \ref{eq2} gives $${\partial\tilde\eta_\theta\over\partial
x_i}(\theta,0)={\partial\tilde\eta_i\over \partial\theta}(\theta,0),$$
$${\partial\tilde\eta_i\over\partial x_j}(\theta,0)={\partial\tilde
\eta_j\over\partial x_i}(\theta,0),$$
we obtain
\begin{eqnarray*}
df_t(\theta,x)&=& \left( \tilde\eta_\theta(\theta,0)+
\sum_i{\partial\tilde\eta_\theta\over \partial
x_i}(\theta,0)x_i\right)d\theta\\
&&+ \sum_i\left(\tilde\eta_i(\theta,0)+\sum_j{\partial
\tilde\eta_i\over \partial x_j}(\theta,0)x_j\right)dx_i\\
&=& \tilde\eta_\theta(\theta,x)d\theta+\sum_i\tilde\eta_i (\theta,x)
dx_i
\end{eqnarray*}
up to first order in $x$.

Damping $f_t$ out to $0$ outside $N(C)$, we arrive at
$\eta_t=\tilde\eta_t-df_t$.  Finally, we obtain the vector field
$X_t$ such that $i_{X_t}\omega_t=\eta_t$.  $X_t$ will then give rise to a
1-parameter family of symplectomorphisms, away from $C$, once we
establish that $X_t\rightarrow 0$ rapidly enough as $p\rightarrow
C$ ($p\in M$).

On $N(C)$, 
\begin{eqnarray}
\omega_t&=& L_1(\theta,x)(d\theta dx_1+dx_2dx_3)\nonumber\\
&+&L_2(\theta,x)(d\theta dx_2+dx_3dx_1)\nonumber\\
&+&L_3(\theta,x)(d\theta dx_3+dx_1dx_2)\nonumber\\
&+& Q, \label{Lij}
\end{eqnarray}
where $L_i(\theta,x)=\sum_jL_{ij}(\theta)x_j$ and $Q$ consists of forms in
$d\theta$ and $dx_i$, whose coefficients are quadratic or higher in the
$x_i$.  In terms of matrices, $\omega_t$ corresponds to 

$$A=\left(
\begin{array}{rrrr}
0 & L_1 & L_2 & L_3\\
-L_1 & 0 & L_3 & -L_2\\
-L_2 & -L_3 & 0 & L_1\\
-L_3 & L_2 & -L_1 & 0\\
\end{array}\right)
+\widetilde Q,$$
where 
$\widetilde{Q}$ has quadratic or higher terms in the $x_i$ and the
matrix is with respect to basis $\{d\theta,dx_1,dx_2,dx_3\}$.
$i_{X_t}\omega_t=\eta_t$ then becomes
$$ 
(a_\theta \mbox{ } a_1 \mbox{ } a_2 \mbox{ } a_3)
A=
(\eta_\theta \mbox{ } \eta_1 \mbox{ } \eta_2 \mbox{ } \eta_3)
$$
with $X_t=a_\theta d\theta +\sum_ia_idx_i$.  Thus,
\begin{eqnarray*}
(a_\theta \mbox{ } a_1 \mbox{ } a_2 \mbox{ } a_3) 
 & = & ( \eta_\theta \mbox{ } \eta_1 \mbox{ } \eta_2 \mbox{ } \eta_3)
 A^{-1}\\
&=& {(\eta_\theta \mbox{ } \eta_1 \mbox{ } \eta_2 \mbox{ } \eta_3)
\over L_1^2+L_2^2+L_3^2}\left(
\begin{array}{rrrr}
0 & -L_1 & -L_2 & -L_3\\
L_1 & 0 & -L_3 & L_2\\
L_2 & L_3 & 0 & -L_1\\
L_3 & -L_2 & L_1 & 0\\
\end{array}\right)
\end{eqnarray*}
up to first order in $x$.  This means that $|X_t|< k|x|$ near $C$;
hence, as $x\rightarrow 0$,
$|\phi_1(\theta,x)-\phi_0(\theta,x)|\rightarrow 0$, where $\phi_t$
is the flow such that ${d\phi_t\over dt}=X_t$. This concludes our
proof. \QED

\section{Local normal forms}

On a neighborhood $N(C)= C\times D^3$ of $C$, $\omega$ can be written as in
Equation \ref{Lij}.  If $\omega$ is generic, then it is transverse to the zero
section of $\ext^+$, and $(L_{ij}(\theta))$ is nondegenerate for all $\theta$.

\begin{lemma}  $(L_{ij}(\theta))$ is symmetric and traceless.
\end{lemma}

\noindent
{\bf Proof:} By comparing 0th order terms in the $x_i$, $d\omega=0$ implies
$$ {\partial L_1\over \partial x_1}+ {\partial L_2\over \partial x_2}
+{\partial L_3\over \partial x_3}=0,$$
$${\partial L_2\over\partial x_3}-{\partial L_3\over \partial x_2}=0, \space\space
{\partial L_3\over \partial x_1}-{\partial L_1\over \partial x_3}=0, \space\space
{\partial L_1\over \partial x_2}-{\partial L_2\over \partial x_1}=0. 
$$
\QED

\vskip.2in
\noindent
$(L_{ij}(\theta))$ thus has a basis $\{v_1(\theta),
v_2(\theta), v_3(\theta)\}$ of eigenvectors  for each $\theta$ (though 
the $v_i$ are not necessarily continuous in
$\theta$).  Since $(L_{ij}(\theta))$ is traceless, either two of the eigenvalues are positive and the remaining is negative for all $\theta$, or vice versa.
Hence, $(L_{ij}(\theta))$ gives rise to a splitting of 
$\R^3\times S^1\rightarrow
S^1$ into a real line bundle over $S^1$ and a rank 2 vector bundle over $S^1$.
Such splittings are classified by homotopy classes of $S^1$ into 
$\R{\mbox{\bf P}}^2$, and $\pi_1(\R{\mbox{\bf P}}^2)=\Z/2$.  Hence,

\begin{prop} There exist two splittings of $\R^3\times S^1\rightarrow S^1$,
the oriented one and the unoriented one.
\end{prop}

What is rather remarkable is that we have the following:

\begin{thm}\label{two}
There exist SD harmonic 2-forms for the product metric on $S^1\times D^3$
for both types of splittings:
\end{thm}

\noindent
(A) 
\vskip-.5in
\begin{eqnarray*}
\omega_A &=& x_1(d\theta dx_1 + dx_2dx_3)\\
&+& x_2(d\theta dx_2 + dx_3dx_1)\\
&-& 2x_3(d\theta dx_3 + dx_1dx_2)\\
&=& *_3\mu + d\theta\wedge \mu,
\end{eqnarray*}
where $\mu=d({1\over 2} (x_1^2+x_2^2)-x_3^2)$, and $*_3$ is the $*$-operator
for the flat metric on $D^3$.  Here, $(L_{ij}(\theta))=\mbox{diag}(1,1,-2)$,
with fixed positive and negative eigenspaces. Note that $\omega_A$ is $S^1$-invariant.

\vskip.15in
\noindent
(B) 
\vskip-.5in
\begin{eqnarray*}
\omega_B &=& (x_1\cos \theta + x_2\sin \theta) e^{x_3}(d\theta dx_1 +dx_2dx_3)\\
& + & (x_1\sin \theta-x_2\cos \theta) e^{x_3}(d\theta dx_2 +dx_3dx_1)\\
& + & R(-x_1(d\theta dx_1 + dx_2dx_3) +x_3(d\theta dx_3 +dx_1dx_2)),
\end{eqnarray*}
with $0<R<1$.  Here,
$$(L_{ij}(\theta))=\left(
\begin{array}{ccc}
\cos \theta-R & \sin \theta & 0\\
\sin \theta& -\cos \theta & 0\\
0 & 0 & R
\end{array}\right).$$
Observe that there are two positive eigenvectors, one of which is $(0,0,1)$.
It is a direct computation to show that $d\omega_B=0$ and that the splitting
of $\R^3$ given by $(L_{ij}(\theta))$ is the {\em unoriented} one.

\vskip.15in
We can alternatively construct an $\omega_B$, which is not the same as the one
above, but arises more naturally.  Starting with 
\begin{eqnarray*}
\omega &=& x_1(d\theta dx_1 + dx_2dx_3)\\
&+& x_2(d\theta dx_2 + dx_3dx_1)\\
&-& 2x_3(d\theta dx_3 + dx_1dx_2)
\end{eqnarray*}
on $[0,2\pi]\times D^3$, we glue $\phi: \{2\pi\}\times D^3 \rightarrow \{0\}\times D^3$ via
\begin{eqnarray*}
\theta & \mapsto  & \theta-2\pi\\
x_1 & \mapsto & x_1\\
x_2 & \mapsto & -x_2\\
x_3 & \mapsto & -x_3.
\end{eqnarray*}
$\phi^*\omega=\omega$, and obtain an $\omega_B$ on $S^1\times S^2$
corresponding to the unoriented splitting.

\begin{thm}\label{local}
 Given an SD harmonic 2-form $\omega$, there exists a 1-parameter
family of perturbations $\{\omega_t\}$, local near $C$,
 such that $\omega_0=\omega$, $\omega_1|_{N(C)}$ is one of the two local forms
as in Theorem \ref{two} (up to $\pm \omega$), and $[\omega_t]\in H^2(M;\R)$ is
independent of $t$.
\end{thm}

\noindent
{\bf Proof:} Given $\omega$, consider the neighborhood $N(C)=C\times D^3$ of
one of the circles.  Assume we are in case (A).  Case (B) is identical.  
After an orthonormal change of frame, we may write
\begin{eqnarray*}
\omega&=& (L_{11}(\theta)x_1+L_{12}(\theta)x_2)(d\theta dx_1 + dx_2dx_3)\\
& + & (L_{21}(\theta)x_1+L_{22}(\theta)x_2)(d\theta dx_2 +dx_3dx_1)\\
& + & \lambda_3(\theta)x_3(d\theta dx_3 + dx_1dx_2)\\
& + &  Q,
\end{eqnarray*}
with, say, $(L_{ij}(\theta))_{1\leq i,j \leq 2}$ positive definite and 
$\lambda_3(\theta)<0$.  Here, the $L_{ij}(\theta)$ and $\lambda_3(\theta)$
are differentiable in $\theta$.

Now, take a 1-parameter family $\omega_t=(1-t)\omega+t\omega_A$ on $N(C)$.  
After shrinking $N(C)$ if necessary, $\omega_t$ is symplectic on $N(C)$ away
from $C$.  Using a local version of Moser's theorem (see the proof of Theorem
\ref{thm:sddiffeo}), we see that there exists a $C^1$-diffeomorphism
$$\phi:(N_0(C),\omega)\stackrel{\sim}{\rightarrow}(N_1(C),\omega_A),$$
where $N_0(C)$, $N_1(C)$ are small neighborhoods of $C$, $\phi=id$ on $C$,
and $\phi$ is a smooth map away from $C$.  Hence $\phi$ allows us to remove
$(N_0(C),\omega)$ and graft on $(N_1(C),\omega_A)$.  We can perform this
operation through a 1-parameter family $\omega_t$, and hence there exists
a global family $\omega_t$ on $M$ with $\omega_0=\omega$ and $\omega_1|_{N(C)}
=\omega_A$.  Moreover, the perturbation can be performed in an arbitrarily 
small neighborhood of $C$ and without altering the cohomology class.\QED

\vskip.2in
In essence, Theorem \ref{local} tells us that, in studying the singular
circles of $\omega$, we may assume that the zeros are either (A) or (B).

\section{Contact structures on the boundaries}

In this section we investigate the boundary properties of $\omega_A$ and 
$\omega_B$. More precisely, we have

\begin{thm}
There exist contact forms $\lambda_A$ and $\lambda_B$ on $\partial N(C)=
S^1\times S^2$ such that $d\lambda_A=\omega_A$ and $d\lambda_B=\omega_B$ on
$S^1\times S^2$.
\end{thm}

\noindent
{\bf Proof:} (A)  For example, consider the following $S^1$-invariant contact
1-form
$$\lambda=-{1\over 2}(x_1^2+x_2^2-2x_3^2)d\theta +x_2x_3dx_1 -x_1x_3dx_2$$
on $N(C)$. We then compute that $d\lambda=\omega$ on $N(C)$ and 
$$\sum_ix_idx_i\wedge\lambda\wedge d\lambda=\left({1\over 2}(x_1^2+x_2^2)(x_1^2
+x_2^2+2x_3^2) + 2x_3^4\right) d\theta dx_1dx_2dx_3.$$
Since $S^1\times S^2=\{\sum_ix_i^2=1\}$ is a leaf of $\sum_ix_idx_i$, 
$i_{S^1\times S^2}^*(\lambda\wedge d\lambda)\not=0$ if and only if 
$\lambda\wedge d\lambda\wedge\sum_ix_idx_i\not=0$ on $S^1\times S^2$.
Here $i_{S^1\times S^2}$ is the inclusion $S^1\times S^2\rightarrow S^1\times
D^3$.  Noting that $\lambda\wedge d\lambda\wedge\sum_ix_idx_i=0$ if and only if
$x_1=x_2=x_3=0$, we have that $i^*_{S^1\times S^2}(\lambda)$ is a contact
1-form on $S^1\times S^2$ with $di^*_{S^1\times S^2}(\lambda)=
i^*_{S^1\times S^2}.$  Thus, $(M-N(S^1),\omega)$ is a symplectic manifold
with contact boundary.

\vskip.2in
\noindent
(B) Consider the 1-form
$$\lambda=-{1\over 2}(x_1^2+x_2^2-2x_3^2)d\theta +x_2x_3dx_1 -x_1x_3dx_2$$
on $[0,2\pi]\times D^3$.  $d\lambda=\omega$, and $\phi^*\lambda=\lambda$, 
where $\phi$ was the glueing map of Theorem \ref{two}, so we glue together 
a contact 1-form $\lambda_B$ such that $d\lambda_B=\omega_B$.  The rest is
the same as (A).
 \QED

\vskip.25in
Let us now describe the orbits of the Reeb vector fields.  
\vskip.15in
\noindent
(A) $\omega_A$ is compatible with a metric $\tilde{g}=\lambda g$, where $g$
is the standard product metric on $S^1\times D^3$.  We can then write the
compatible $J$ satisfying $\tilde{g} (x,y)=\omega(Jx,y)$ as $J={1\over \lambda}
A$, where 
$$A=\left(\begin{array}{cccc}
0 & x_1& x_2 & -2x_3 \\
-x_1 & 0 & -2x_3 & -x_2\\
-x_2 & 2x_3 & 0 & x_1\\
2x_3 & x_2 & -x_1 & 0
\end{array}\right)$$
represents $\omega$ with respect to $\{\theta, x_1,x_2,x_3\}$. Now, the Reeb
vector field $X$ for $i^*_{S^1\times S^2}\omega$ is given, up to multiple,
by 
$$J(\sum_ix_i{\partial\over \partial x_i})={1\over \sqrt{x_1^2+x_2^2+4x_3^2}}
(x_1^2+x_2^2-2x_3^2, -3x_2x_3, 3x_1x_3, 0)^T.$$

Finally, $\lambda_A(X)=1$ implies that
$$ X={1\over f}\left[ (x_1^2+x_2^2-2x_3^2){\partial\over \partial \theta}
-3x_2x_3{\partial\over \partial x_1}+3x_1x_3{\partial \over \partial x_2}\right],$$
with $$f=-{1\over 2}\left[(x_1^2+x_2^2)(x_1^2+x_2^2+2x_3^2)+2x_3^4\right].$$

Solving for the orbits, $x_1^2+x_2^2$ and $x_3^2$ are fixed for each orbit, and
hence,
\begin{eqnarray*}
x_1 &=& \sqrt{1-r^2} \cos{R_1(r) t}\\
x_2 &=& \sqrt{1-r^2} \sin{R_1(r) t}\\
x_3 &=& r\\
\theta &=& R_2(r) t +c,
\end{eqnarray*}
where $r$ is a constant, and $R_1$ and $R_2$ are functions of $r$.  

In particular, the noteworthy closed orbits are 
$\{(0,0,1)\}\times S^1$, $\{(0,0,-1)\}\times S^1$, and $\{(x_1,x_2,0)\}\times 
S^1$, with $x_1^2+x_2^2=1$ and $x_1$, $x_2$ fixed.  These correspond to the
stable and unstable gradient directions in the Morse theory of ${1\over 2}
(x_1^2+x_2^2-2x_3^2)$ near $(0,0,0)$.  Moreover, the orbit $\{(0,0,1)\}\times 
S^1$ is {\em nondegenerate}, and so is the family $\{(x_1,x_2,0)\}\times
S^1$.  There are other closed orbits, but these do not seem to have any
Morse-theoretic significance.

\vskip.15in
\noindent
(B) We apply the previous considerations and work on $[0,2\pi]\times S^1/\sim$.
There is one orbit $\{(0,0,\pm 1)\}\times S^1$, which is a double of the 
orbits for (A).    Since
$\phi$ identifies $(2\pi, (x_1,x_2,0)) \sim (0, (x_1,-x_2,0))$, we also have the doubled 
 closed orbits $\{(x_1, \pm x_2,0)\}\times S^1$, with $x_2\not = 0$, and the 
single closed orbits $\{(1,0,0)\}\times S^1$, $\{(-1,0,0)\}\times S^1$.  

\vskip.2in
\noindent
{\bf Remark:} There is an example of a singularity of type (B) bounding
a disk, which can be made to vanish.

\noindent
\hskip.3in\small{Mathematics Department, Princeton University, Princeton, NJ 08544}

\noindent
\hskip.3in\small{{\em honda@math.princeton.edu}}
\end{document}